I: INTRODUCTION

In a previous paper (Gottlob, 2004) some inconsistencies in the formalism of the modern Bayes' Theorem (BT) were noted. These and other inconsistencies triggered further investigations leading to a new formalism that will be derived in this paper. BT was regarded to be the only method available for increasing the probability of a hypothesis due to some new evidence (C. Glymour, 1980, p. 91). The same author conjectured that "there must be relations between evidence and hypotheses that are important to scientific argument and to confirmation but to which the Bayesian scheme has not yet penetrated". J. Earman, author of an encompassing treatise on Bayes' Theorem (Bayes or Bust?, 1992) confessed on p. l that he was a Bayesian on Mondays, Wednesdays and Fridays, however, that the rest of the week, he had doubts about the viability of the system as the basis for analyzing scientific inference. Nevertheless, I am aware of the fact that BT was generally in use for about seven decades and that but a tiny minority of the large number of users uttered objections. Therefore, I start right here with an example that will show that with the Bayesian formalism not everything can be in perfect order:

According to a rule that is considered to be an important part of BT, if a hypothesis (H), implies some evidence, (E), the probability $P(H|E)$ will be $P(H) / P(E)$. However, from this follows the absurd conclusion that if the evidence is certain, $P(H|E)$ will be $P(H) / 1 = P(H)$. In other words: Certain evidence will not be able to confirm a hypothesis to any degree. Almost certain evidence [$P(E)$ close to 1], will contribute but minimally to confirmation.

Here, a practical *Example:* A doctor sees a patient suffering from high fever. He assumes 'malaria', however, other febrile states cannot be excluded. The doctor assigns his hypothesis (H) the probability .5. He decides to examine several blood smears and assigns a positive result (plasmodia in blood cells) probability 1. And in fact, he does find plasmodia, germs that are causing the disease, but are not found in healthy persons. One might object that in BT, E indicates the *antecedent* probability of evidence, a prediction that is experimentally confirmed later. However, the argument of time dependence will not help in our example. The doctor assigns the evidence `plasmodia in the blood´ probability 1 for malaria *before* examination of the blood smear. Should the probability of malaria after detecting plasmodia be .5 / 1 = .5?[1]

---

[1] This example cannot be disconfirmed by C. Howson's (1984, 1985) reply to similar objections, e. g., by C. Glymour (1980) and by R. Campbell and T. Vinci (1983). Howson claims that in BT, P(e) should not be relativized to the whole background knowledge K, but only to K - {e} – However, this proposal is at odds with the generally accepted rule that in probabilistic reasoning our total background knowledge must be considered. See also D. Christensen, 1999

Because of the above mentioned absurdity, a completely new scheme, the rMPE-Method, mainly based on three principles so far not sufficiently considered by Bayesianism, was designed. These principles are: (i): The cMPE-Method for additions of probabilities in spite of their non-linearity (Gottlob, 2000, 2002, 2003); (ii): A consequent distinction between weight-bearing and extensional evidence (Gottlob. 2004), and (iii): avoiding increasing probabilities by divisions by variables that distort the result. With this method, we can assign any hypothesis the probability according to the contemporary *actual* evidence. If the actual evidence is 1, this will cause an increase of P(H) to P(H) $^+$ P(H), i.e., P(H) added to P(H), however allowing for the non-linearity of the system. The criticism of BT will focus on

(i)  The absurdity outlined in Example 1.
(ii)  The lacking differentiation between *'weight-bearing evidence ($E_W$)'* and *' extensional evidence'* ($E_{EXT}$) in BT. (See also Gottlob, 2004).
(iii)  The necessity in BT to increase probabilities by *divisions,* whereby the divisor is a variable < 1. This procedure causes additional non-linearities to probability systems that are non-linear by nature and causes diminished probability results. Divisions by constant divisors or multiplications do not yield these drawbacks.

Of lesser or undefined importance may be:

(iv) The derivation of BT from the conjunctive axiom that in turn is derived from a *descending* formula. i.e., the result of the operation is lower than the single premises, whereas due to mutual support the result should increase.
(v)  The impossibility (or at least difficulty) to check the correctness of probabilistic hypotheses or algorithms by internal proofs or by proofs according to reality.
(vi) Tacit conditions that are not reflected by the BT formalism.

*Some notations and definitions:*

P(H): The probability of hypothesis H. P(H|E): The probability of H given E.
As evidence E, the evidence is considered that actively supports or inhibits some probabilities, e.g, of hypothesis H, whereby the evidence for the occurrence of E may be relevant, too. P(A) $^+$ P(B); P(A) $^-$ P(B): Addition, respectively, subtraction of the probability of A by the probability of B, however, taking into account non-linearity. The probabilities of error, 1- P(A) , 1- P(B), the complements of P(A) and P(B), are also notated as P(¬A) and P(¬B) or P($A^C$) and P($B^C$). — All probabilistic reasoning must allow for our background knowledge. For sake of clarity this is dropped in our formulas, however it is always considered tacitly. In Ch. II the cMPE-Method is briefly described. Ch. III will be devoted to the criticism of BT. In Ch. IV, the rMPE-Method will be derived and presented.
In Ch. V, a general discussion and the description of some special applications will follow.

II: THE cMPE-METHOD

was derived for increasing probabilities by addition, however allowing for their nonlinearity. As the method is already described (Gottlob, 2000, 2002, 2003 and in arXiv 2004) I shall mention here only what is necessary for understanding the principle.

*a) The non-linearity of probabilities*

Probabilities can range only between 0 and 1. Thus simply adding probabilities, e. g., .4 +.7 cannot be carried out. Various probabilities may belong to quite different sample spaces containing, say, the inhabitants of a village and of a whole country. In this case, no simple addition is feasible, even if the sum would be below 1. There must also be a greater difference between probability .9 and certainty than between probabilities .8 and .9. Fig. 1 will show these differences (see also Fig. 4 in Gottlob, 2004). - Only a minority of probability operations works in a linear way. The precondition therefore is: Belonging to one sample space and, accordingly, real or possible sums not exceeding 1. These conditions are met by *disjunctions*.

*b) The Principle of Complementarity.*
By multiplications of numbers <1 we arrive at products that are smaller than the factors. Besides simple addition there are two possibilities for arriving at results that exceed the premises: *(i) Divisions by numbers < 1*. As will be shown in Ch. IV, Fig. 5, such divisions cause non-linearities, even if the dividend remains constant and the divisor increases linearly. This is different in multiplications that, in the bargain, have the advantage of commutability.
(ii) An alternative would be the *multiplication of the complements* of various probabilities. Hereby, the product will be smaller than the factors; however, by subtracting the product from 1, the difference will be larger than any of the original probabilities. By this multiplicative procedure, no additional non-linearity is caused (see Fig. 1, below).
Thus we may add two or more probabilities by multiplying their complements and by subtracting the product from 1. A precondition is that the probabilities support one another. One more precondition is semantic independence: All added probabilities must differ by the way of observation, (e.g., by different sense organs), or by being based on different scientific methods, by being characterized by different properties, by being observed by independent witnesses, or being supported by different circumstantial facts. Only the conclusion, they support must be common (See Fig. 7, below, for a model of semantic independence). Of course, if not the whole probability supports the other, only the effective part of that probability may enter into the calculation.
On the other hand, there might be other influences that inhibit the support. Here, the most important may be pieces of evidence E, that are not exclusively related to our hypothesis but are found anywhere else, too: $P(E|\neg H)$

*c) Derivation of the MPE-Method.*
The probability of the occurrence of two or more completely independent events by sheer chance is
$$P(A \ \& \ B \ \& \ \ldots \ \& \ N) = P(A) \times P(B) \times \ldots \times P(N) \tag{1}$$
This is the special multiplication rule. The general multiplication rule for two probabilities reads:
$$P(A \ \& \ B) = P(A) \times P(B|A) = P(A|B) \times P(B) \tag{2}$$
(2) is valid for conditional probabilities, e.g., $P(A|B)$, where $P(A) \neq P(A|B)$. Both, (1) and (2) are *descending* formulas, because due to the coincidence, the resulting chance is smaller than the initial values. Therefore, these formulas cannot be applied, if A, B or more events support (corroborate) one another. In this case, the final probability must be greater than the probability of the single event, and (1) or (2) will be of no use. Only by the cMPE-Method we may add probabilities, taking into account their non-linearity and avoiding sums >1. For computing `ascending' probabilities, we may first multiply the probabilities of error, $P(A^C)$ and $P(B^C)$. This, again, is a descending operation. The product of $P(A^C) \times P(B^C)$ or of more probabilities of error becomes smaller than the premises. The decrease in probability is caused by multiplying numbers < 1.
$$P(\neg A \ \& \ \neg B \ \& \ldots \& \ \neg N) = P(A^C) \times P(B^C) \times \ldots \times P(N^C) \tag{3}$$
I named (3), the 'MPE-Method' (**M**ultiplication of the **P**robabilities of **E**rror). This method may play an important role in our everyday cognitions. If we see, smell, and taste an apple, all sensations support one another and the probability of error becomes smaller, almost vanishing. As already mentioned, a basic requirement is that P(A), P(B) etc. are semantically independent.

*d: Derivation of the cMPE- Method.*

If we subtract the result of the MPE- Method from 1, we arrive at the common probability of our cognition or of the occurrence of some events that support one another. For indicating an ascending calculation, we notate the support between the items by a superscript $^+$.

$$P(A)\ ^+ P(B) = 1 - P(A^C) \times P(B^C) \qquad (4a)$$

For more evidential facts, supporting one another, and of course, if no counterexample is known, we have:

$$P(A)\ ^+ P(B)\ ^+ .... \ ^+ P(N) = 1 - P(A^C) \times P(B^C) \times .... \times P(N^C) \qquad (4b)$$

If we carry out the multiplication in (4a), we arrive at (4c):

$$P(A)\ ^+ P(B) = P(A) + P(B) - P(A) \times P(B) \qquad (4c)$$

By the MPE-Method probabilities are added, however, allowing for their non-linearity[2]. Here, for a proof, I shall limit myself to Fig. 1:

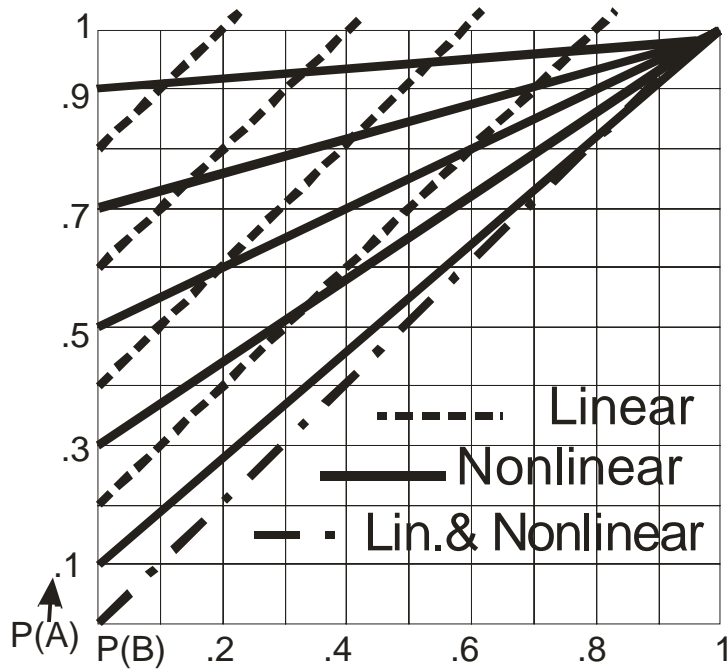

Fig. 1: The data of the abscissa are added to the probability, indicated on the ordinate. Continuous lines: cMPE-values. Broken lines: Linear additions of probabilities.

---

[2] For further proof, the reader is referred to Gottlob, (2002, 2003), where also the difference of (4c) to the union and to the probability that at least one event of a multitude of events must occur, is discussed. In 2003, also the conversion of (4a – c) into percentages is considered.

In Fig 1 we see that by linear adding of the probabilities indicated on the abscissa to the probabilities indicated on the ordinate (broken lines) straight parallel lines are generated that may intersect the horizontal 1 line at different points. Adding probabilities non-linearly to one of the probabilities indicated on the ordinate generates straight lines that converge at 1.0. The non-linearity is indicated by different slopes.

Transformations of cMPE-operations into percentage operations were treated in 2003, where also a percentage formula of P.O. Ekelöf´s (1964) and a severe criticism by L.J. Cohen (1976, 1977) is discussed.

The graph Fig. 1 may be regarded as proof for a non-linear addition, because the lines start at the point P(B) = 0, their course must be straight, because of one factor remaining constant, and must end at 1, where all lines must converge, even if probability 1 is added. If one term of the sum, say P(B) is assumed to be zero, the result will equal the other term, P(A): P(A) + 0 = P(A).

In Table 1, P.8 will be added by the MPE-Method to various probabilities.

Table 1
.8 is added to the Initial value.

| Initial Value | Sum | initial Value | Sum |
|---|---|---|---|
| ,1 | .82 | .9 | .98 |
| .3 | .86 | .99 | .998 |
| .5 | .9 | .999 | .999 8 |
| .8 | .96 | .999 9 | .999 98 |

We see in Tab. 1 that the impact of adding .8 by a nonlinear method to various probabilities decreases strongly with increasing initial value and that the resulting probabilities approach 1 asymptotically. Probability 1 can be reached only, if 1.0 is added.

For completeness, in (5), the formula for non-linear *subtraction* (DPE-Method, **D**ivision of the **P**robabilities of **E**rror) will be stated, whereby the superscript again indicates non-linearity:

$$P(A) \bar{\ } P(B) \bar{\ } P(C) \bar{\ } ... \bar{\ } P(N) = 1 - P(A)^C / P(B^c) \times P(C^C) \times .... \times P(N^c) \qquad (5)$$

Now we may compute the difference between .999 98 and .999 8. Linearly we arrive at .000 18, however, by the DPE-Method we arrive at 1-.00002 /.0002 = .9, the 5000fold value.

Because of this non-linearity, we must doubt, that the support function, S, frequently mentioned in literature, e.g. "S(H,E) = P(H|E) - P(H)" can be correct. Here, the DPE-Method should be applied: P(H|E) $\bar{\ }$ P(E) will reveal the real difference.

Fig. 2 will show the ascending operation (b and c) as opposed to the descending operation (a).

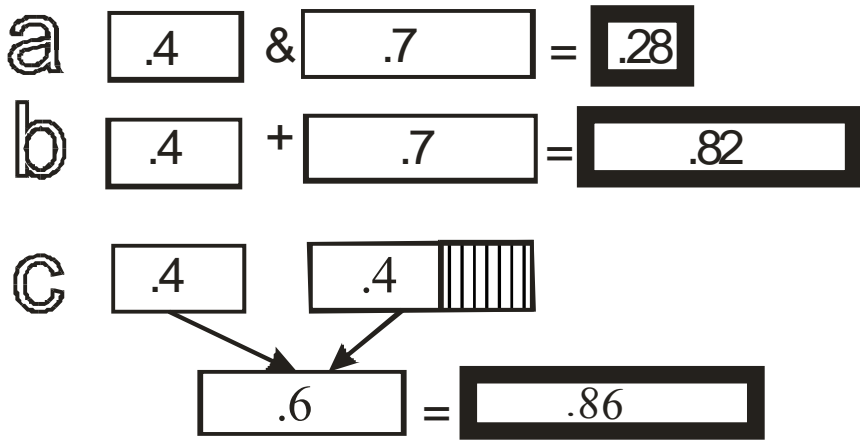

Fig. 2. a: Premises, P(A) =.4, P(B) =.7, P(A & B) = .28.
  b: cMPE- Method, $.4 \,{}^+ .7 = 1 - ,6 \times .3 = .82$
  c: A, totally (P = .4) and B partially, with P = .4 support C (P = .6);
    The bold faced quadrangles indicate the result.

III: THE NEO-BAYESIAN THEOREM AND ITS CRITICISM.
*a) The Theorem.*
Originally (6a) was deduced from (2), the conjunctive axiom. According to Russell (1948, p. 364), from this follows the 'principle of inverse probability'(6a):
P(H|E) = P(H) x P(E|H) / P(E)     (6a)
Frequently, especially in papers of earlier decades, (6a) is also called Bayes' Theorem. In this equation an increasing P(E) reduces the probability, however, if P(E) is 1, P(H|E) becomes P(H) . Another fundamental drawback of this formula consists in its not allowing for the probability of evidence, found under the condition of ¬H. This drawback could be remedied by considering the total probability in the denominator. We shall discuss here the most simple case that P(H|E) is just calculated by the prior: P(H), the likelihood P(E|H), P(E) and the probability of P(E|¬H). In this case (6b) reads:
P(H|E) = P(H) x P(EH) / [P(H) x P(E|H) + P(¬H) x P(E|¬H)]     (6b)
If H implies E, P(E|H) becomes 1 and (6a) becomes (6c):
    P(H|E) = P(H) / P(E)     (6c)

*b) Criticism*
(i) The absurdity that in (6c), *P(H/E) cannot be increased, if P(E) =1*, was already mentioned.
(ii) An important point is the *lacking differentiation between extensional evidence ($E_{EXT}$) and weight-bearing evidence ($E_W$)*[3] in (6a) and in (6c). These different concepts of evidence were dealt with in Ch. V of Gottlob (2004). Any piece of evidence that bears on the confirmation of hypothesis H may consist of two constituents: P($E_W$), the part that increases P(H) and P(E)¬P(H)), the part that reduces P(H) by being found not related to the hypothesis. Both of the parts jointly form P($E_{EXT}$), frequently just notated `E´.

---

[3] The term `weight-bearing evidence´ was coined in order to differ from older concepts of weight of evidence, proposed by J.M. Keynes (1957, Ch.5) or by I.J. Good, 1950, Ch.. 6.

*(iii)* The necessity in BT, *to increase probabilities by division.* True, by divisions by numbers <1, the quotient increases. However, there might arise some inconsistencies by the natural differences between divisions and multiplication. Unlike in multiplications, the constituents of divisions are non-commutative and the quotients, obtained from linear series of divisors are non-linear. Example: Dividend: .8; divisors; .2, .4, .6. .8. The respective quotients are 4, 2, 1.33, 1. In addition, the positive role of P($E_W$) and the negative role of P(E) – P(H) in increasing P(H) to P(H|E) is not specified. We shall come back to the problem of non-linear quotients in Ch. IV, Fig. 5.

*(iv)* (6a - c) are *derived from the conjunctive axiom* (2). This axiom, also called the 'general multiplication rule' indicates a descending operation, P(H) x P(E|H) or P(H|E) x P(E) are both smaller than the factors and this decrease is due to the multiplication of numbers < 1. The decrease is understandable, because it expresses changes by chance. However, for P(H|E) with few exceptions, we want to compute the *increase of* P(H) due to *"active support",* and therefore we apply a completely different algorithm, based on the ascending cMPE-operation.

(v) In contradistinction to the rMPE-method, to be described below, a convincing proof of the result of the Bayesian operations is impossible. True: in an extreme case, if P(H) implies P(E) and both probabilities are equal, P(H|E) by laws of mathematical logic must become 1, and this is correctly the case for (6a, b. c). However, at the other extreme, [P(E = 1)], we have the difficulties described. Checking the results by experiments in reality is theoretically possible but extremely difficult and checking 'from within' is not feasible. In contradistinction thereto, the cMPE - and the rMPE-Method may be checked by a graphical proof, in a way, a proof 'from within'.

(vi) There are *tacit conditions that are not reflected by the formalism*, e.g., concerning background knowledge 'K', such as not allowing "relativizing the probabilities in the support definition to the entire set K[.... ] including [evidence] e, whereas they should be relativized to K - {e}" (C. Howson, 1985). And the stipulation that "the Bayesian assesses the contemporary support E gives H by how much the agent would change his odds on *H were he now* [counterfactually] to come to know E" (C. Howson, 1984) [brackets added by R.G].

These here criticized features of Bayesianism can be overcome by the rMPE-Method described in Ch. IV.

IV: CONSTRUCTING AN ALTERNATIVE ALGORITHM BASED ON THE COMPLEMENTARITY PRINCIPLE AND ON WEIGHT-BEARING EVIDENCE, HOWEVER, AVOIDING DIVISIONS.

The inventors of the modern BT had to resort to divisions, because by multiplications of probabilities, the result will be reduced instead of being increased and no method for adding probabilities, in spite of their non-linearity, was known. However, the 'principle of complementarity' described above, permits to add probabilities by multiplying their complements and subtracting the sum from 1, according to (4 a - c). For sake of simplicity, we start with.

*a) Implication of E by H.*

The derivation of the rMPE-Method is demonstrated by Figs. 3 – 5.

Our basic assumption is that P(H|E), in the conventional understanding, may be computed by P(H) $^+$ P($E_W$). As Fig. 4 shows, P($E_W$) = 1 – [P($E_{EXT}$) – P)H)]:

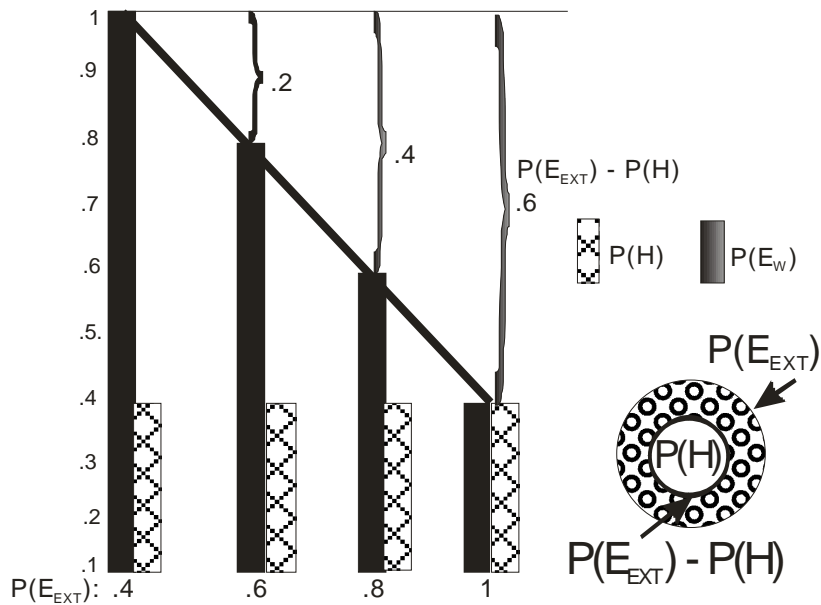

Fig. 3 Derivation of the rMPE-Method for implications of E by H. For hypothesis H in all instances, we assume .4. On the right a Venn diagram, showing the implication.

The following formulas serve for the derivation of the rMPE-formula:
<a> $P(E_{EXT})$ & $\neg H) = P(E_{EXT}) - P(H)$
See the small Venn-Diagram on the right side of Fig. 3. A linear subtraction is permissible, because, due to the implication, the data belong to the same sample space.
<b> $P(E_W) = 1 - [P(E_{EXT}) - P(H)]$
See Fig. 3.
<c> Only $P(E_{EXT}) - P(H)$ has some negative influence on $P(H)$, it reduces $P(E_W)$,
i.e. the complement of $P(E_{EXT}) - P(H)$ (see <a>).
We now have a look at the various $P(E_W)$s in Fig 3.
If $P(E)$ and $P(H)$ are congruent, $P(H|E)$, as well, as $P(E|H)$ are necessarily 1. But by the cMPE-Method a sum can be only 1, if one term added is 1. Thus, $P(E_W)$ must be 1. This proves (d) of the following considerations:
<d> $P(E_{EXT}) = .4$; $P(E_{EXT}) - P(H) = 0$; $P(E_W) = 1 - 0 = 1$ and $P(H) ^+ P(E_W) = 1$
<e> $P(E_{EXT}) = .6$; $P(E_{EXT}) - P(H) = .2$; $P(E_W) = 1 - .2 = .8$ and $P(H) ^+ P(E_W) = .88$
<f> $P(E_{EXT}) = .8$; $P(E_{EXT}) - P(H) = .4$; $P(E_W) = 1 - .4 = .6$ and $P(H) ^+ P(E_W) = .76$
<g> $P(E_{EXT}) = 1$; $P(E_{EXT}) - P(H) = .6$; $P(E_W) = 1 - .6 = .4$ and $P(H) ^+ P(E_W) = .64$

Also <g> is so by necessity: If $P(E_{EXT}) = 1$, $P(E) - P(H)$ must be $1 - P(H)$ <a> and $P(E_W)$ $1 - 1 - P(H) = P(H)$ <b>). Therefore, $P(H) ^+ P(E_W) = 1 - P(H^C) \times [1 - [P(E_{EXT}) - P(H)]^C$.
All of the four results <d> – <g>) represent $P(H) ^+ P(E_W)$. They form a linear sequence.
From these considerations, we may construct the formula (7) for $P(H) ^+ P(E_W)$:

$P(H)\ ^+ P(E_W) = 1- P(H^C) \times \{1 – [P(E_{EXT}) –P(H)]\}^C = 1 – P(H^C) \times [P(E_{EXT})]$ (7)

For comparison, we compute P(H|E) according to BT (6c): We arrive for P(H) / P(E) at 1, .66, .5 and .4. which is a non-linear sequence, explaining, why evidence 1.0 frequently cannot confirm P(H), as already noted in the introduction. We may conclude from these findings:

*1.) That P(H|E) [computed if P(H) implies P(E) by P(H) / P(E)] is not the correct operation for calculating the increase of P(H) due to evidence P(E) and*

*2.) That the correct operation is $P(H)\ ^+ P(E_W)$. In what follows, this parameter will be notated `$P(H_W)$´, in order to permit a distinction between the rMPE Method and BT.*

$P(H_W)$ can be attained by cMPE-addition of P(H) to $P(E_W)$. Because the latter is also 1- [$P(E_{EXT})$ - P(H)] (formula b, above), we arrive at $P(H_W) = 1- P(H^C) \times 1- 1- [P(E) – P(H)]^C$ so that in (7), P(H) increases and P(E) – P(H) decreases the result, however, the latter factor exerts also an increasing effect, because $P(E_{EXT}) – P(H)$ is 1- $P(E_W)$ and this is contained in (7) as complement, so that it also contributes positively to P(H|E).

(7) is valid for implications of E by H, according to B. Russell (1948 p.364) )"The most useful case". In these examples, p(E) equals or exceeds P(H).

The important advantage of (7), as compared with the BT-implication formula (6c) is that also P(E) = 1 increases the probability of ($H_W$). Indeed, if we assume that in Fig. 3 and 4 the impact of E on P(H) [1- P(E|¬H)] =.4, then, according to (7), $P(H_w)$ must be $.4\ ^+.4 = 1- .6 \times .6 = .64$. In Fig. 4, three typical constellations are compared:

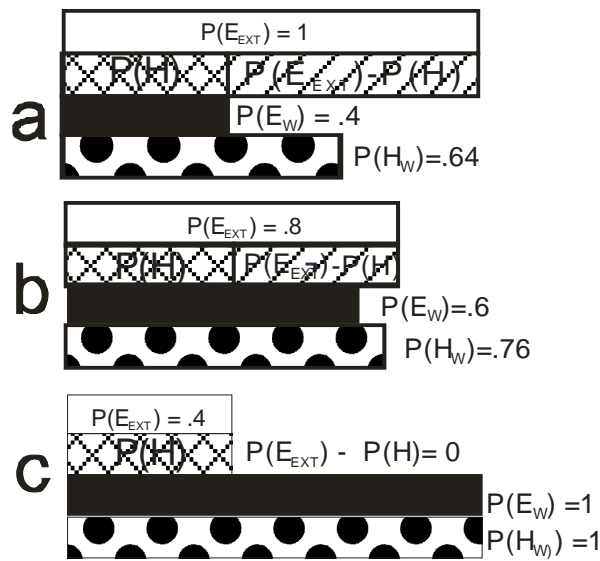

Fig. 4: P(H) = .4. The black bars indicate $P(E_W)$. **a**: $P(E_{EXT}) = 1$, $P(E_W) = P(H)$. **b**: $P(E_{EXT}) = .8$, $P(E_W) = .6$; **c**: $P(E_{EXT}) = .4$, $P(E_W)$ as well as $P(H_W) = 1$.

That (7) is not only valid for P(H) = .4, may be seen on Fig. 5, where a multitude of rMPE-calculations is plotted and also the difference between the rMPE-Method and the conventional formula for implications (6c) is demonstrated.

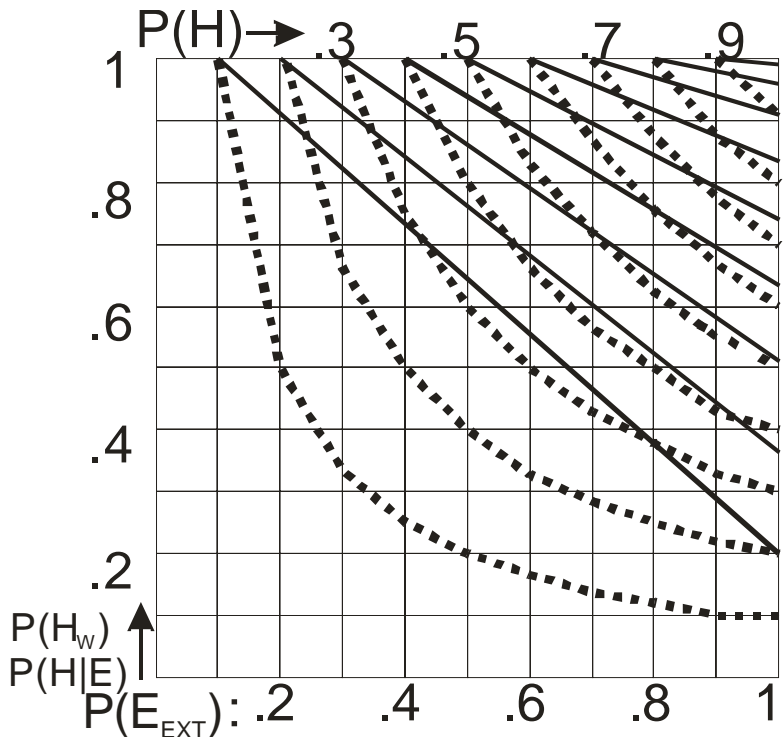

Fig:5: Implication of E by H. The dotted lines indicate P(H) / P(E) according to (6c). The continuous lines indicate $P(H_W)$ according to (7). All lines have their origin at P(H) at the top.

Examples: (The conventional results of BT for implication in cornered brackets):
$P(H) = .3$, $P(E_{EXT}) = .4$. $P(H_W) = 1 - .7 \times .1 = .93$; [P(H) / P(E) = .75]
$P(H) = .6$, $P(E_{EXT}) = .8$. $P(H_W) = 1 - .4 \times .2 = .92$;. [P(H) / P(E) = .75]
All values of $P(H_W)$, with P(E) = 1, result in $P(H)^+ P(H)$, e.g.: P(H) =.5, $P(E_{EXT}) = 1$,
$P(H_W) = 1 - .5 \times (1 - .5) = 1 - .5 \times .5 = .75$.
Again, we may use the results achieved in Fig. 5 as a graphic proof of the method. The initial values of $P(H_W)$ are 1 (as proved in Fig. 3) because there P(H) and P(E) are equal. That the end points equal P(H) + P(H), was explained by Figs. 3 and 4; the straightness of the lines follows from one factor in the multiplication being kept constant.
In contrast to the computation according to (7), the broken lines in Fig. 5, resulting from the division of P(H) by $P(E_{EXT})$, (6c) are curved and non-linear, in spite of the dividend in each curve being kept constant. The non-linearity is caused by different divisors. Please note that for implications, all results of BT are considerably lower than the results of the rMPE-Method for the same initial P(H) value. This is true also for higher values: If P(H) =.7 and $p(E_{EXT}) =1$, we obtain by BT .7 and by the rMPEMethod a $P(H_W)$ of .91, a remarkable difference.

*b: H. implies several pieces of evidence that beyond H do not overlap.*
In 2004, I mentioned some inconsistencies in the procedure put forward by C.D. Broad (1918). The latter procedure consists in dividing probability H by various evidential data. The denominator, according to Broad, reads:
$c,c_2.....cn|f' = c_1|f \times c_2|c_1f \times c_3|c2c1f \times .... c_n| c_{n-1}......c_i|f$. In this denominator $c_1, c_2$ are propositions that are implied by a hypothesis. 'f' indicates background knowledge.
As already remarked by J. Earman (1992, p. 107), any additional c decreases the denominator, which eventually may become smaller than the numerator, and probabilities >1 may be obtained. In the first version of my 2004 paper, I relied upon divisions for increasing probabilities and on the cMPE-Method for adding probabilistic data. As shown in this paper, because of drawbacks with divisions by variables, I now recommend applying the rMPE-Method for these cases. Fig. 6 will explain the new procedure.

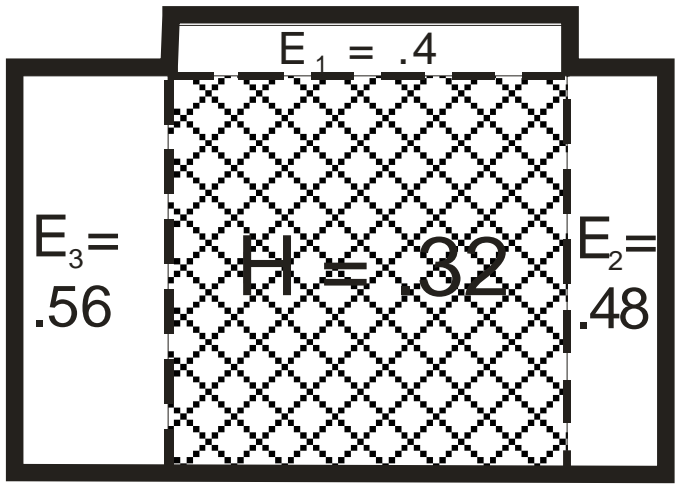

Fig. 6: $P(H) = .32$ is supported by $E_1 = .40$, $E_2 = .48$, and $E_3 = .56$. The $P(E_{EXT}) - P(H)$ are respective .08, .16 and .24. The $P(E_{123})$ are independent beyond H, however completely superimposed on H. -The figure illustrates the principle of semantic independence.

Due to the superimposition, the area of H is covered by three more layers of $E_1$, $E_2$, $E_3$. On the other hand, there are three different $P(E) - P(H)$. Thus the formula for a multitude of semantically independent pieces of evidence, supporting H reads:

$P(H_W | E_1, E_2, \ldots E_n) = 1 - P(H^C) \times [P(E_1)-P(H)] \times [P(E_2)-P(H)] \times \ldots \times [P(E_n)-P(H)]$  (8)

And for the data contained in Fig. 7, we compute: $1 - P(H_W)^C \times [P(E_{EXT1})-p(H)] \times [P(E_{EXT2})-P(H)] \times , [P(E_{EXT3}) - P(H)] = .997\ 9$. At first blush, .997 9 appears rather high. However, we must consider that P(H) is supported by three semantically independent probabilities and that the $P(E) - P(H)$ thereof are reduced by mutual multiplication..

*c: Conditional probabilities without implication.*
Here the standard BT formula (6b) seems to produce correct results, because the same results are arrived at also by quite simple procedures, such as the relative frequency formula (9)
        P = Number of favorable cases / Number of possible cases         (9).
I show this by an example of R.L. Winkler(1972, p.43):
A person, chosen at random from a certain population is given a tuberculin skin test.
In persons having tuberculosis, the skin test will be positive in 98%: P(pos|disease) = .98.
Persons, not suffering from tuberculosis have a positive skin test in 5%: P(false pos) = .05
1 % of the population has tuberculosis. What is the probability that a person with a positive skin test has tuberculosis? From the premises we may calculate:

| | | |
|---|---|---|
| P(H x P(E\|H) | .01 x .98= | .0098 |
| 99 persons without the disease, but 0.5 with positive reaction: | .99 x .05 = | .0495 |
| Sum of positive skin tests in 100 persons: | | .0593 |
| Quotient: (P disease and positive test) / ( P of all positive reactions): | .0098 / .0593 = | .165 |

R.L. Winkler arrives at the same result by BT (6b).The calculative effort is almost equal. However, the advantage of (6b), over (6a) or (6c) is that if $P(E) - PH$ is zero, $P(H|E)$ becomes 1.This advantage (6b) shares with (7) and (9). On the other hand, cases with $P(E) = 1$ (or very close to 1) can be correctly computed only with (7).

Does the conformity of the results of (6b) with the somewhat more trivial relative frequency calculation prove the correctness of (6b)? Not necessarily, because both of the methods suffer from the same drawback, they are based on the same division that causes aggravation of the non- linearity, similarly, as depicted in Fig. 5. A method would be desirable that permits to add the probability of the hypothesis non-linearly to the weight-bearing capacity of the evidence also in non-implication cases.

V: DISCUSSION AND SPECIAL CASES.
*a: The following flaws or inconsistencies in BT could be avoided:*
(i) If H implies E, according to (6c), a $P(E_{EXT}) = 1$, no longer fails to support the hypothesis, but confirms it according to $P(H_W) = P(H) \,^+ P(H)$.
(ii) There is clear differentiation between the concepts of $P(E_W)$ and of $P(E_{EXT})$, concepts that are to a great part antagonistic.
(iii) Divisions by variable divisors that aggravate the natural non-linearity of probabilities are avoided.
(iv) The method is based on ascending operations (as it should be), viz.: on the cMPE Method and not based on descending procedures, such as the conjunctive axiom. Descending probabilities as results of multiplications, may be converted into ascending ones by embedding them into complements and subtracting the result from 1.
(v) The formulas derived for the rMPE-Method can be proved graphically (Fig. 3) by drawing straight lines from a starting point at $P(H_W) = 1$, where $P(H_{EXT}) = P(E)$, to a logically determined end point, $[P(H) \,^+ P(H)]$, whereby the straightness of the lines is given by necessity. This 'internal proof' is not feasible with BT.
(vi) Besides semantic independence, no special tacit rules not supported by the formalism must be observed. There is no time factor and no counterfactual assumption prescribed. Any knowledge available at the time of calculation may be used for computing the final probability.

*b: Application in practical and theoretical science.*
The rMPE-Method may be used in all cases that until now were tackled by BT. This concerns conditional probability cases and implication cases that until now were solved by (6c). In combination with the cMPE-Method, in all fields of science common probabilities of semantically independent data may be evaluated. This may be useful in practical science as well as in Theory of Science. By the cMPE-Method, two or more synergistic probabilistic data may be added and hereby evaluated without any likelihood (P(E|H)) known (See Gottlob, 2004 for more examples). Here an example by D. Christensen (1999), however, in my interpretation:

An agent is wondering whether deer is in a nearby wood. He comes across a pile of deer droppings and considers the probability of error for `deer in the wood´ very low, say, .001.

Shortly thereafter, he finds a shed antler. Again, the probability of error is very low (nobody was interested to put the antler to that place, in order to deceive him). The probability of error is at most again .001. By the MPE-Method he arrives at .001 x .001 = .000 001 and by the cMPE-Method at at least .999 999. The increase in probability due to the second finding is only .000 999, however, due to the non-linearity, this correlates with an increase by .999. Of course, these calculations are not performed numerically. Rather tacitly or even unconsciously, the agent reasons: Two semantically independent observations with very small probability of error support one another. From this results a vanishing probability of error and, by order of magnitude, the highly certain assumption that there is deer in the wood. This example shows that by the cMPE-Method in practice probabilities may be increased, without knowing any likelihood. In addition, for implications of E by H, the rMPE-Method may be applied for all questions that until now were tried to be solved by BT. Besides many advantages, the rMPE-Method described here, is easily manageable and the probabilities arrived at, frequently will be considerably higher than the results, obtained with BT.

*c:.• Astonishing or surprising evidence.*

A classic example: B. Russell (1948, p. 365) mentions the confirmation of the law of gravitation. There were seven planets known. However, because of some anomalies in the orbit of Uranus, an eighth planet was assumed and its localization predicted. Indeed, Neptune was found exactly at the site determined beforehand. How can these findings be accommodated to our formalism? Astronomers subdivide the heaven by rectascensions and declinations into 41 253 square areas. Thus, the chance of finding - without the theory of gravitation - another planet in a certain area, $P(E|\neg H)$, was 7/ 41 253 = .000 17, and the already very high probability of the theory of gravitation, say, .999 9, becomes according to (4a) 1-.000 1 x .000 17 = .999 999 98 or $\approx 1-10^{-8}$. -These and other considerations favor the assumption that the psychological moments of surprise or astonishment, because hardly numerically recordable, are of less importance than the improbability of finding the facts in question. In other words: Astonishing things are astonishing, because of the small $P(E|\neg H)$, the small chance of encountering them, unless H were true.

*d:- The old evidence problem.*

A possible solution of this problem by the cMPE-Method and by the DPE-Method was briefly outlined by R. Gottlob (2004). This question will be dealt with in extension elsewhere.

*e: Multiple evidence.*

According to the rMPE-Method (8) and as shown in Fig.6, the probabilities of error of multiple pieces of evidence that are semantically independent are multiplied in the course of processing and thereby are vigorously reduced. Any new factor will cause additional reduction. After subtraction of the product from 1, the resulting probabilities may be high and frequently, no longer distinguishable from 1 for all practical purposes.

ACKNOWLEDGEMENT: The author feels highly indebted for many enlightening discussions to Prof. Viktor Scheiber, former head of the Institute of Medical Statistics, Vienna University.
Author's E-Mail: rainer.gottlob@univie.ac.at